# IRRATIONALITY MEASURES, IRRATIONALITY BASES, AND A THEOREM OF JARNÍK

JONATHAN SONDOW

ABSTRACT. We recall that the irrationality exponent $\mu(\alpha)$ of an irrational number $\alpha$ is defined using the irrationality measure $1/q^\mu$. Using $1/\beta^q$ instead, and motivated by a conditional result on Euler's constant, we define the *irrationality base* $\beta(\alpha)$. We give formulas for $\mu(\alpha)$ and $\beta(\alpha)$ in terms of the continued fraction expansion of $\alpha$. The formulas yield bounds on $\mu(\alpha)$ and $\beta(\alpha)$ involving the Fibonacci numbers, as well as a growth condition on the partial quotients of $\alpha$ implying $\mu(\alpha)=2$, and a weaker condition implying $\beta(\alpha)=1$. A theorem of Jarník on Diophantine approximation leads to numbers with prescribed irrationality measure. By a different method, we explicitly construct series with prescribed irrationality base. Many examples are given.

## 1. INTRODUCTION

Recall that the irrationality exponent $\mu(\alpha)$ of an irrational number $\alpha$ is defined in terms of the irrationality measure $1/q^\mu$. Using $1/\beta^q$ instead, we introduce a weaker measure of irrationality, the *irrationality base* $\beta(\alpha)$, as follows. If there exists a real number $\beta \geq 1$ with the property that for any $\varepsilon > 0$ there is a positive integer $q(\varepsilon)$ such that

$$\left|\alpha - \frac{p}{q}\right| > \frac{1}{(\beta+\varepsilon)^q},$$

for all integers $p, q$, with $q \geq q(\varepsilon)$, then the least such $\beta$ is called the irrationality base $\beta(\alpha)$ of $\alpha$.

The motivation for using the irrationality measure $1/\beta^q$, rather than some other function of $q$, is the following. In [**6**] we gave criteria for irrationality of Euler's constant, $\gamma$, involving the Beukers-type double integral

$$I_n : \iint_{[0,1]^2} -\frac{(x(1-x)y(1-y))^n}{(1-xy)\log xy}\,dx\,dy = \binom{2n}{n}\gamma + L_n - A_n,$$

where

$$d_{2n}L_n \in \mathbf{Z}\ln(n+1) + \mathbf{Z}\ln(n+2) + \cdots + \mathbf{Z}\ln(2n)$$

and $d_{2n}A_n \in \mathbf{Z}$, with $d_n$ denoting the least common multiple of the numbers $1, 2, \ldots, n$. In particular, we showed that $\gamma$ is irrational if the fractional part of $d_{2n}L_n$ exceeds $2^{-n}$ infinitely often. In a paper in preparation, we prove that a stronger condition implies a measure of irrationality for $\gamma$. Namely, if

(C) $$\lim_{n \to \infty} \frac{\ln \| d_{2n} L_n \|}{n} = 0,$$

where $\| t \|$ denotes the distance from $t$ to the nearest integer, then the estimates as $n \to \infty$

$$I_n = 4^{-2n(1+o(1))}, \quad \binom{2n}{n} = 2^{2n(1+o(1))}, \quad d_{2n} = e^{2n(1+o(1))}$$

imply that for any $\varepsilon > 0$ there exists $q(\varepsilon) > 0$ such that

$$\left| \gamma - \frac{p}{q} \right| > \frac{1}{(2e + \varepsilon)^q},$$

for all integers $p, q$, with $q \geq q(\varepsilon)$; in other words, the irrationality base of $\gamma$ satisfies the inequality $\beta(\gamma) \leq 2e = 5.436\ldots$. We also show that a condition similar to (C), but based on an integral for $\ln(\pi/2)$ [**7**], implies that $\ln \pi$ is irrational and $\beta(\ln \pi) \leq 2$. Moreover, we present numerical evidence for both conditions. These results suggest that it would be interesting to study the new irrationality measure $1/\beta^q$ by itself.

In Section 2 we define irrationality measures, exponents, and bases. In Section 3 we prove formulas for $\mu(\alpha)$ and $\beta(\alpha)$ in terms of the continued fraction expansion of $\alpha$, and compute several examples. The formulas yield bounds on $\mu(\alpha)$ and $\beta(\alpha)$ involving the Fibonacci numbers, as well as a growth condition on the partial quotients of $\alpha$ implying $\mu(\alpha) = 2$, and a weaker condition implying $\beta(\alpha) = 1$. In Section 4 we observe that the existence of numbers with prescribed irrationality measure follows from V. Jarník's 1931 construction of continued fractions with certain Diophantine approximation properties [**3**]. In the final section, we use a different method to construct explicitly series with prescribed irrationality base. The two methods yield numbers with the same irrationality base but different approximation properties.

## 2. IRRATIONALITY MEASURES, EXPONENTS, AND BASES

After recalling the definitions of irrationality measure and irrationality exponent, we introduce the notion of irrationality base. Throughout this section and the next, $\alpha$ denotes a fixed but arbitrary irrational number.

**Definition 1.** An *irrationality measure* is a function $f(x, \lambda)$, defined for $x \geq 1$ and $\lambda > 0$, which takes values in the positive reals and is strictly decreasing in both $x$ and $\lambda$. If there exists $\lambda > 0$ with the property that for any $\varepsilon > 0$ there exists a positive integer $q(\varepsilon)$ such that

$$\left|\alpha - \frac{p}{q}\right| > f(q, \lambda + \varepsilon), \quad \text{for all integers } p, q, \text{ with } q \geq q(\varepsilon),$$

then we denote by $\lambda(\alpha)$ the least such $\lambda$, and we say that $\alpha$ *has irrationality measure* $f(x, \lambda(\alpha))$. Otherwise, if no such $\lambda$ exists, we write $\lambda(\alpha) = \infty$.

**Definition 2.** Let $f(x, \lambda) = x^{-\lambda}$. If $\alpha$ has irrationality measure $f(x, \mu) = x^{-\mu}$ for some $\mu = \mu(\alpha)$, then $\mu(\alpha) \in [2, \infty)$ is called the *irrationality exponent* of $\alpha$. If there is no such $\mu$, then we write $\mu(\alpha) = \infty$ and we call $\alpha$ a *Liouville number.*

*Remark.* Our definition of irrationality exponent is equivalent to the one in [**2**, p. 106], but our definition of irrationality measure differs from that in [**2**, p. 17].

**Definition 3.** Let $f(x, \lambda) = \lambda^{-x}$. If $\alpha$ has irrationality measure $f(x, \beta) = \beta^{-x}$ for some $\beta = \beta(\alpha)$ (so that $\beta$ is the smallest number with the property that for any $\varepsilon > 0$ there exists $q(\varepsilon) > 0$ such that

(1) $$\left|\alpha - \frac{p}{q}\right| > \frac{1}{(\beta + \varepsilon)^q}, \quad \text{for all integers } p, q, \text{ with } q \geq q(\varepsilon)),$$

then we call $\beta(\alpha) \in [1, \infty)$ the *irrationality base* of $\alpha$. Otherwise, if no such $\beta$ exists, we write $\beta(\alpha) = \infty$ and we say that $\alpha$ is a *super Liouville number.*

**Example 1.** *The sum of the following series is a super Liouville number*:

$$S = \frac{1}{1} + \frac{1}{2^1} + \frac{1}{4^{2^1}} + \cdots.$$

*Proof.* If we write $S = \sum_{n \geq 1} b_n^{-1}$, where $b_1 = 1$ and $b_n = (2^n)^{b_{n-1}}$, for $n > 1$, then the $n$th partial sum of the series equals $a_n / b_n$, for some integer $a_n$, and we have

$$0 < S - \frac{a_n}{b_n} < \frac{2}{b_{n+1}} = \frac{2}{(2^{n+1})^{b_n}} < \frac{1}{(2^n)^{b_n}} < \frac{1}{b_n^2}, \quad n \geq 2.$$

By the last bound, $S$ is irrational. By the next-to-last bound, no number $\beta$ can satisfy the requirements of Definition 3 with $\alpha = S$. So $\beta(S) = \infty$. ●

Writing inequality (1) as

$$\left|\alpha - \frac{p}{q}\right| > q^{-q \log(\beta + \varepsilon)/\log q}$$

in order to compare it with the irrationality measure $q^{-\mu}$ in Definition 2, we see that an upper bound on the irrationality base of a number is a much weaker condition than an upper bound on its irrationality exponent. In particular, we have the following proposition, which justifies the terminology in the last part of Definition 3.

**Proposition 1.** *An irrational number which is not a Liouville number has irrationality base one; equivalently, if $\beta(\alpha) > 1$, then $\mu(\alpha) = \infty$. In particular, a super Liouville number is also a Liouville number.*

*Proof.* If $\mu = \mu(\alpha)$ is finite, then for any $\varepsilon > 0$ we have

$$\left| \alpha - \frac{p}{q} \right| > \frac{1}{q^{\mu+\varepsilon}} > \frac{1}{(1+\varepsilon)^q},$$

for all integers $p, q$, with $q > 0$ sufficiently large, and so $\beta(\alpha) = 1$. In particular, $\beta(\alpha) = \infty$ implies $\mu(\alpha) = \infty$. ●

**Example 2.** *The transcendental numbers $\ln 2$, $e$, $\pi$, $\zeta(2)$, and the irrational number $\zeta(3)$, all have irrationality base one.*

*Proof.* Upper bounds on their irrationality exponents are known [**2**, Chapter 2, §§4, 5]. ●

**Corollary 1.** *An algebraic irrational number has irrationality base one.*

*Proof.* This is a immediate from Liouville's theorem. ●

### 3. FORMULAS FOR THE IRRATIONALITY EXPONENT AND BASE

We prove formulas for the irrationality exponent and base of an irrational number $\alpha$ in terms of its continued fraction expansion. Then we give two corollaries and several examples. The first corollary gives bounds on $\mu(\alpha)$ and $\beta(\alpha)$ involving the Fibonacci numbers. The second gives two growth conditions on the partial quotients of $\alpha$: one implies that $\mu(\alpha) = 2$, and the other, weaker one that $\beta(\alpha) = 1$.

In this section and the next, we will use the following two lemmas.

**Lemma 1.** (Legendre [**5**]) *For integers $p, q$, with $q > 0$, the inequality*

$$\left| \alpha - \frac{p}{q} \right| < \frac{1}{2q^2}$$

*implies that $p/q$ is a convergent of the continued fraction expansion of $\alpha$.*

*Proof.* See [**4**, p. 11]. ●

**Lemma 2.** *For any real numbers $C > 0$, $\lambda > 1$ there exists a positive integer $q_0 = q_0(C, \lambda)$ such that if $p, q$ are integers satisfying the inequalities*

$$\left| \alpha - \frac{p}{q} \right| \leq \frac{C}{\lambda^q}, \quad q \geq q_0,$$

*then $p/q$ is a convergent of $\alpha$.*

*Proof.* By Lemma 1, it suffices to take $q_0$ so large that $C/\lambda^{q_0} < 1/(2q_0^2)$. ●

Denote a simple continued fraction by

$$[b_0; b_1, b_2, \ldots] = b_0 + \cfrac{1}{b_1 + \cfrac{1}{b_2 + \cdots}}.$$

**Theorem 1.** *The irrationality exponent and irrationality base of an irrational number $\alpha$ with continued fraction expansion $\alpha = [b_0; b_1, b_2, \ldots]$ and convergents $p_n/q_n$ are given by*

(2) $$\mu(\alpha) = 1 + \limsup_{n \to \infty} \frac{\ln q_{n+1}}{\ln q_n} = 2 + \limsup_{n \to \infty} \frac{\ln b_{n+1}}{\ln q_n},$$

(3) $$\ln \beta(\alpha) = \limsup_{n \to \infty} \frac{\ln q_{n+1}}{q_n} = \limsup_{n \to \infty} \frac{\ln b_{n+1}}{q_n}.$$

*Remark.* The formulas in (2) for the irrationality exponent are probably known, but we have not found a reference.

*Proof.* If we define $\lambda_n$ by the equation

(4) $$\left| \alpha - \frac{p_n}{q_n} \right| = q_n^{-\lambda_n},$$

then the second of the inequalities [**4**, p. 8]

(5) $$\frac{1}{2 q_n q_{n+1}} < \left| \alpha - \frac{p_n}{q_n} \right| < \frac{1}{q_n q_{n+1}}$$

implies that $\lambda_n > 2$. Using (5) and Lemma 1, we deduce that

$$\limsup_{n \to \infty} \lambda_n = \mu(\alpha).$$

Taking logarithms in (4) and (5) leads to the first formula in (2), and the second follows using the relations [**4**, p. 2]

(6) $$q_{n+1} = b_{n+1} q_n + q_{n-1} = b_{n+1} q_n (1 + o(1)).$$

To prove (3), we define $\lambda_n$ by replacing the right-hand side of (4) by $\lambda_n^{-q_n}$. Then (5) and Lemma 2 imply that $\limsup_{n \to \infty} \lambda_n = \beta(\alpha)$, and (3) follows as did (2). ●

**Example 3.** *The golden mean $\phi = (1 + \sqrt{5})/2$ has irrationality exponent* 2 *and irrationality base* 1.

*Proof.* Here $q_n = F_n$ is asymptotic to $\phi^n$, and so $\mu(\phi) = 1 + 1 = 2$ and $\ln \beta(\phi) = 0$. ●

**Corollary 2.** *If $\alpha = [b_0; b_1, b_2, \ldots]$ and $F_n$ is the n-th Fibonacci number, then*

(7) $$\limsup_{n \to \infty} \frac{\ln b_{n+1}}{\ln(F_n b_1 b_2 \cdots b_n)} \leq \mu(\alpha) - 2 \leq \limsup_{n \to \infty} \frac{\ln b_{n+1}}{\ln(F_n + b_1 b_2 \cdots b_n)},$$

(8) $$\limsup_{n \to \infty} \frac{\ln b_{n+1}}{F_n b_1 b_2 \cdots b_n} \leq \ln \beta(\alpha) \leq \limsup_{n \to \infty} \frac{\ln b_{n+1}}{F_n + b_1 b_2 \cdots b_n}.$$

*Proof.* Using $q_1 = b_1$, the recursion in (6), and induction (or the fact that $q_n$ is the sum of $F_n$ products of the $b_i$), we obtain the bounds

$$F_n + b_1 b_2 \cdots b_n - 1 \leq q_n \leq F_n b_1 b_2 \cdots b_n,$$

and (7), (8) follow, using (2), (3). ●

**Example 4.** *The numbers*

$$L_1 = [0; 2^{1!}, 2^{2!}, 2^{3!}, \ldots], \quad L_2 = [0; 1, 2^1, 4^{2^1}, \ldots],$$

$$L_3 = [0; 2, 2^2, 2^{2^2}, \ldots], \quad L_4 = [0; 1, 2^1, 3^{2^1}, \ldots]$$

*all have irrationality exponent infinity and irrationality base one. In particular, they are Liouville numbers, but not super Liouville numbers. On the other hand, the numbers*

$$S_1 = [0; 1, 2^2, 3^{3^3}, \ldots], \quad S_2 = [0; 1, 2^2, 4^{4^4}, \ldots]$$

*are super Liouville numbers.*

*Proof.* We weaken the lower bound in (7) and the upper bound in (8) to

$$\mu(\alpha) \geq \limsup_{n \to \infty} \frac{\ln b_{n+1}}{\ln(2^n b_1 b_2 \cdots b_n)}, \quad \ln \beta(\alpha) \leq \limsup_{n \to \infty} \frac{\ln b_{n+1}}{b_1 b_2 \cdots b_n}.$$

Then

$$\mu(L_1) \geq \limsup_{n \to \infty} \frac{(n+1)!}{n + 1! + 2! + \cdots + n!} \geq \lim_{n \to \infty} \frac{(n+1)!}{3n!} = \infty,$$

$$\ln \beta(L_1) \leq \limsup_{n \to \infty} \frac{(n+1)! \ln 2}{2^{1!+2!+\cdots+n!}} \leq \lim_{n \to \infty} \frac{(n+1)!}{2^{n!}} = 0.$$

For $L_2$ we have $b_n = (2^n)^{b_{n-1}}$, and so

$$\mu(L_2) \geq \limsup_{n \to \infty} \frac{(n+1)b_n}{n + 1 + 2b_1 + 3b_2 + \cdots + nb_{n-1}} \geq \lim_{n \to \infty} \frac{2^{b_{n-1}}}{1 + nb_{n-1}} = \infty,$$

$$\ln \beta(L_2) \leq \limsup_{n \to \infty} \frac{(n+1)\ln 2}{b_1 b_2 \cdots b_{n-1}} \leq \lim_{n \to \infty} \frac{n+1}{2^{n-2}} = 0.$$

The estimates for $L_3$ and $L_4$ are similar. For $S_1$, use the tower notation $T_n(k)$, defined by

(a) $$T_1(k) = k, \quad T_{n+1}(k) = k^{T_n(k)},$$

to write $b_n = T_n(n)$. Since

$$2^n T_1(1) T_2(2) \cdots T_n(n) < n^{T_{n-1}(n+1)}, \quad n \geq 3,$$

by (8) we have

$$\ln \beta(S_1) \geq \limsup_{n \to \infty} \frac{T_n(n+1) \ln(n+1)}{2^n T_1(1) T_2(2) \cdots T_n(n)}$$

$$\geq \lim_{n \to \infty} \frac{T_n(n+1)}{n^{T_{n-1}(n+1)}} = \lim_{n \to \infty} \left(\frac{n+1}{n}\right)^{T_{n-1}(n+1)} = \infty.$$

The estimate for $S_2$, where $b_n = T_n(2^{n-1})$, is easier. ●

*Remark.* Any one of the numbers $L_1, L_2, L_3, L_4$ shows that the converse of Proposition 1 is false.

**Corollary 3.** *If the partial quotients $b_0, b_1, b_2, \ldots$ of $\alpha$ satisfy the growth restriction $b_n = e^{o(n)}$ as $n \to \infty$ (for example, $b_n$ bounded, or a polynomial in n, or $b_n < e^{n/\ln n}$), then $\mu(\alpha) = 2$. If $b_n = e^{o(F_n)}$ as $n \to \infty$ (for instance, $b_n < e^{(3/2)^n}$), then $\beta(\alpha) = 1$.*

*Proof.* By (7), (8), we have

$$0 \leq \mu(\alpha) - 2 \leq \limsup_{n \to \infty} \frac{\ln b_{n+1}}{\ln F_n}, \quad 0 \leq \ln \beta(\alpha) \leq \limsup_{n \to \infty} \frac{\ln b_{n+1}}{F_n}.$$

Since $\ln F_n = n \ln \phi - O(1)$ as $n \to \infty$ (see Example 3), the first statement follows. The second statement is immediate (note that $3/2 < \phi$). ●

## 4. PRESCRIBING IRRATIONALITY MEASURES VIA JARNÍK'S THEOREM

In 1931, V. Jarník proved a general result on simultaneous Diophantine approximation of $s \geq 1$ numbers [**3**, Satz 6] (see also [**1**, p. 18, Exercise 1.5]). We recall the statement for $s = 1$ and give his construction by continued fractions (Theorem 2). We then point out that Jarník's method produces numbers with prescribed irrationality measure (Theorem 3 and Corollary 4).

In the next section, bypassing Jarník's method, we explicitly construct series with irrationality base any prescribed real number $\beta > 1$. As a preview of the two methods, here are numbers with the same irrationality base but different approximation properties, one a continued fraction and the other a series.

**Example 5.** *We have $\beta(\theta_2) = 2$ for the continued fraction $\theta_2 = [0; 3, 1, 2, 17, b_5, b_6, \ldots]$, where $b_{n+1}$ is related to the n-th convergent $p_n/q_n$ by the formula*

$$(9) \qquad b_{n+1} = 1 + \left\lfloor \frac{2^{q_n}}{q_n^2} \right\rfloor.$$

*In particular, $b_5 = 1 + \lfloor 2^{191}/191^2 \rfloor > 10^{52}$.*

*Proof.* Taking $\omega(x) = 2^{-x}$ in Corollary 4 below, formula (9) satisfies conditions (10). ●

**Example 6.** *We have $\beta(\tau_2) = 2$ for the series*

$$\tau_2 = \frac{1}{2} + \frac{1}{2^2} + \frac{1}{2^{2^2}} + \cdots .$$

*Proof.* Take $\beta = 2$ in (12) below. ●

*Remark.* Although the numbers in Examples 5 and 6 have the same irrationality base $\beta(\theta_2) = \beta(\tau_2) = 2$, they have opposite approximation properties: for any integers $p, q$, with $q > 0$ sufficiently large, we will see that

$$\left|\theta_2 - \frac{p}{q}\right| < \frac{1}{2^q}, \quad \left|\tau_2 - \frac{p}{q}\right| > \frac{1}{2^q}.$$

In order to state Jarník's theorem, we need a definition. (We use the terminology of [**1**, Chapter I, §3].)

**Definition 4.** (Jarník) Given a real number $\theta$ and a real-valued function $\omega(x)$, defined and positive for $x \geq 1$, we say that $\theta$ *is approximable to order* $\omega$ if for every $C > 0$ there exists a pair of integers $p, q$, with $q > C$, such that

$$\left|\theta - \frac{p}{q}\right| < \omega(q).$$

**Theorem 2.** (Jarník) *Suppose that the function $\omega(x)$ is defined for $x \geq 1$, positive, and decreasing, and that $\omega(x) = o(x^{-2})$ as $x \to \infty$. Then there exists an irrational number $\theta$ which is approximable to order $\omega$, but not to order $c\omega$, for any constant $0 < c < 1$.*

*Proof.* (Sketch) Jarník constructs a continued fraction by choosing the partial quotients $b_1, b_2, \ldots$ successively so that $b_{n+1}$ is related to the $n$th convergent $p_n/q_n$ by the formulas

(10) $$b_{n+1} \omega(q_n) q_n^2 > 1, \quad \lim_{n \to \infty} b_{n+1} \omega(q_n) q_n^2 = 1.$$

Then he shows that any number $\theta = [0; b_1, b_2, \ldots]$ which satisfies (10) also satisfies the conclusions of the theorem. ●

The following is a simple application of Jarník's theorem.

**Theorem 3**. *Suppose that $f(x, \lambda)$ is an irrationality measure with a fixed value of $\lambda > 0$ such that $f(x, \lambda) = o(x^{-2})$ as $x \to \infty$. Suppose further that for any $\varepsilon > 0$ there exists a positive constant $c = c(\varepsilon) < 1$, independent of $x$, such that $f(x, \lambda + \varepsilon) < cf(x, \lambda)$, for $x$ sufficiently large. Then there exists an irrational number $\theta_\lambda$ which has irrationality measure $f(x, \lambda)$ and is approximable to order $\omega(x) = f(x, \lambda)$.*

*Remark.* Although an irrationality measure is strictly decreasing in both variables, the second hypothesis is not redundant, because $c$ does not depend on $x$.

*Proof.* The function $\omega(x) = f(x,\lambda)$ satisfies the hypotheses of Theorem 2. Let $\theta_\lambda$ be a number satisfying the conclusions. Since $f(x,\lambda)$, and thus $\omega(x)$, is a decreasing function of $x$, and $\theta_\lambda$ approximable to order $\omega$, it follows that $\lambda(\theta_\lambda) \geq \lambda$ (in the notation of Definition 1). To show that $\lambda(\theta_\lambda) \leq \lambda$, suppose on the contrary that $\lambda(\theta_\lambda) > \lambda$. Then there exist $\varepsilon > 0$ and infinitely many integers $p, q$, with $q > 0$, such that

$$\left| \theta_\lambda - \frac{p}{q} \right| \leq f(q, \lambda + \varepsilon) < cf(q,\lambda) = c\omega(q),$$

for some constant $0 < c < 1$. But then $\theta_\lambda$ is approximable to order $c\omega$, which is a contradiction. Hence $\lambda(\theta_\lambda) = \lambda$ and the result follows. ●

**Corollary 4.** *Every real number $\beta > 1$ is the irrationality base of some number $\theta_\beta$, which is approximable to order $\beta^{-x}$; any continued fraction $\theta_\beta = [0; b_1, b_2, \ldots]$ satisfying* (10) *with $\omega(x) = \beta^{-x}$ will do. The analogous result for $\mu > 2$ and the irrationality exponent holds with $\omega(x) = x^{-\mu}$.*

*Proof.* This follows from Theorem 3 and the proof of Theorem 2. ●

## 5. PRESCRIBING IRRATIONALITY BASES VIA SERIES

By a different method, we prove a related existence result for the irrationality base.

**Theorem 4.** *Any real number $\beta > 1$ can be realized as the irrationality base of an explicitly constructed series $\tau_\beta$ which is* not *approximable to order $\beta^{-x}$; in fact,*

(11) $$\left| \tau_\beta - \frac{p}{q} \right| > \frac{1}{\beta^q}, \quad \text{for all integers } p, q, \text{ with } q > 0 \text{ sufficiently large.}$$

*If $\beta = b/a > 1$ is a rational number in lowest terms, with $b > 0$, then*

(12) $$\tau_\beta = \beta^{-1} + \beta^{-b} + \beta^{-b^b} + \cdots.$$

To prove Theorem 4, we first prove a lemma in which we construct a series with prescribed irrationality base from a sequence of rational numbers with certain properties. We then prove the theorem by constructing the required sequence. The constructions per se do not involve continued fractions, but the proof of the lemma does.

**Lemma 3.** *Given a real number $\beta > 1$, suppose that*

$$r_n = \frac{a_n}{b_n}, \quad n \geq 1,$$

*is a sequence of reduced fractions with positive denominators satisfying the conditions*

(13) $$1 > r_1 \geq r_2 \geq r_3 \geq \cdots \geq \beta^{-1},$$

(14) $$\lim_{n \to \infty} r_n = \beta^{-1},$$

(15) $$b_n \mid b_{n+1}, \quad n \geq 1, \quad \text{and}$$

(16) $$b_{n+1} \leq t_n, \quad n \geq 2,$$

*where $t_n$ is the tower of height n defined by*

(17) $$t_0 = 1, \quad t_n = b_n^{t_{n-1}} = b_n^{b_{n-1}^{\cdot^{\cdot^{\cdot^{b_1}}}}}, \quad n \geq 1.$$

*Then the sum of the series*

(18) $$\tau_\beta = \sum_{n=1}^{\infty} r_n^{t_{n-1}} = \frac{a_1}{b_1} + \left(\frac{a_2}{b_2}\right)^{b_1} + \left(\frac{a_3}{b_3}\right)^{b_2^{b_1}} + \cdots$$

*has irrationality base $\beta$. Moreover, $\tau_\beta$ satisfies (11).*

*Proof.* Conditions (13) and (17) imply that the series converges. Since $(a_n, b_n) = 1$, for $n \geq 1$, conditions (15) and (17) imply that the $n$th partial sum of the series has denominator $t_n$, so that for some integer $u_n$,

(19) $$s_n := \sum_{i=1}^{n} r_i^{t_{i-1}} = \sum_{i=1}^{n} \frac{a_i^{t_{i-1}}}{t_i} = \frac{u_n}{t_n}, \quad (u_n, t_n) = 1.$$

From (13), (17), and (18), it follows that there is a constant $C > 1$ such that

(20) $$r_{n+1}^{t_n} < \tau_\beta - s_n < C r_{n+1}^{t_n}, \quad n \geq 1.$$

Using (19), (20), and (14), we deduce that if $0 < \varepsilon < \beta - 1$, then

$$(21) \qquad 0 < \tau_\beta - \frac{u_n}{t_n} < \frac{1}{(\beta-\varepsilon)^{t_n}} < \frac{1}{t_n^2}, \qquad n \geq n(\varepsilon).$$

It follows that $\tau_\beta$ is irrational, and that $\beta(\tau_\beta) \geq \beta$. To show that $\beta(\tau_\beta) \leq \beta$, we prove the stronger condition (11). By Lemma 2, it suffices to show that the convergents $p_m/q_m$ of $\tau_\beta$ satisfy (11).

By (19), (21), and Lemma 2, all but a finite number of partial sums are convergents,

$$(22) \qquad s_n = \frac{p_{m_n}}{q_{m_n}}, \qquad n \geq n_0,$$

where

$$(23) \qquad q_{m_n} = t_n.$$

Using (23), (22), (20), and (13), we have

$$\frac{1}{2t_n q_{m_n+1}} = \frac{1}{2q_{m_n} q_{m_n+1}} < \tau_\beta - \frac{p_{m_n}}{q_{m_n}} < C r_1^{t_n}, \qquad n \geq n_0.$$

Since $0 < r_1 < 1$, there exists $n_1 \geq n_0$ such that if $n \geq n_1$ then

$$q_{m_n+1} > \frac{r_1^{-t_n}}{2C t_n} > \frac{\log\left(2 t_n^{2t_n}\right)}{\log \beta}.$$

Now consider a convergent $p_m/q_m$ of $\tau_\beta$, with $m$ large. If $p_m/q_m$ is a partial sum (22), then (13) and the lower bound in (20) imply that $p_m/q_m$ satisfies inequality (11). If $p_m/q_m$ is not a partial sum, then from (22) we have $m_n + 1 \leq m < m_{n+1}$, for some $n = n(m) \geq n_1$. Using (16) and (23), we obtain

$$\frac{\log\left(2 b_{n+1}^{2t_n}\right)}{\log \beta} < q_{m_n+1} \leq q_m < q_{m_{n+1}} = t_{n+1}.$$

Hence for $m$, and thus $n$, sufficiently large,

$$\left|\tau_\beta - \frac{p_m}{q_m}\right| \geq \left|\frac{p_{m_{n+1}}}{q_{m_{n+1}}} - \frac{p_m}{q_m}\right| - \left|\tau_\beta - \frac{p_{m_{n+1}}}{q_{m_{n+1}}}\right| > \frac{1}{q_{m_{n+1}} q_m} - C r_1^{t_{n+1}}$$

$$> \frac{1}{t_{n+1}^2} - \frac{1}{2 t_{n+1}^2} = \frac{1}{2 b_{n+1}^{2 t_n}} > \frac{1}{\beta^{q_m}},$$

so that $p_m/q_m$ satisfies inequality (11). This completes the proof of the lemma. ●

*Proof of Theorem 4.* Given $\beta > 1$, construct a sequence of rational numbers $r_n$, as follows.

*Case* 1. If $\beta$ is rational, say $\beta = b/a$, with $0 < a < b$ and $(a, b) = 1$, take the constant sequence $r_n = a/b$, for $n \geq 1$.

*Case* 2. If $\beta$ is irrational, express $1 - \beta^{-1}$ in "base 2" as

$$1 - \frac{1}{\beta} = \sum_{n=1}^{\infty} \frac{1}{2^{k_n}},$$

with $0 < k_1 < k_2 < k_3 < \cdots$. Write the $n$th partial sum of the series as $S_n = A_n/B_n$, where $A_n$ is odd and $B_n = 2^{k_n}$. Choose a positive integer $n_1$ such that $B_2 \leq T_{n_1}$, where $T_n$ denotes the tower of height $n$ defined by

$$T_0 = 1, \quad T_n = B_1^{T_{n-1}}, \quad 1 \leq n \leq n_1;$$

for $1 \leq n \leq n_1$, set $r_n = 1 - S_1$. Now choose $n_2 > n_1$ such that $B_3 \leq T_{n_2}$, where

$$T_n = B_2^{T_{n-1}}, \quad n_1 + 1 \leq n \leq n_2;$$

for $n_1 + 1 \leq n \leq n_2$, set $r_n = 1 - S_2$. Continuing in this way, we define $r_n$ for all $n \geq 1$.

In both Cases 1 and 2, the sequence $r_1, r_2, r_3, \ldots$ satisfies conditions (13) to (17) (in Case 2, for (17) we have $t_n = T_n$). The theorem now follows from Lemma 3. ●

**ACKNOWLEDGEMENTS.** I am grateful to Yann Bugeaud, Yuri Nesterenko, Tanguy Rivoal, Michel Waldschmidt, and Wadim Zudilin for valuable suggestions and encouragement.

209 WEST 97TH STREET, NEW YORK, NY 11025, USA
*email address*: jsondow@alumni.princeton.edu